\documentclass[pdflatex,sn-mathphys-num]{sn-jnl}

\usepackage{graphicx}%
\usepackage{multirow}%
\usepackage{amsmath,amssymb,amsfonts}%
\usepackage{amsthm}%
\usepackage{mathrsfs}%
\usepackage[title]{appendix}%
\usepackage{xcolor}%
\usepackage{textcomp}%
\usepackage{manyfoot}%
\usepackage{booktabs}%
\usepackage{algorithm}%
\usepackage{algorithmicx}%
\usepackage{algpseudocode}%
\usepackage{listings}%

\usepackage{tikz}
\usepackage{comment}
\usepackage{enumitem}
\setlist[itemize,enumerate]{itemsep=0.5em}
\usepackage{hyperref}
\theoremstyle{plain}
\newtheorem{theorem}{Theorem}

\newtheorem{corollary}[theorem]{Corollary}

\theoremstyle{definition}
\newtheorem{definition}[theorem]{Definition}

\newtheorem*{question}{Question}

\newtheorem{conjecture}{Conjecture}
\newtheorem{proposition}[theorem]{Proposition}

\theoremstyle{remark}
\newtheorem*{remark}{Remark}

\long\def\rev#1{{\color{black}#1}}

\begin{document}
\title{3-manifold polynomials}

\author[1]{\fnm{José} \sur{Frías}}\email{frias@ciencias.unam.mx}

\author[2]{\fnm{José Carlos} \sur{Gómez-Larrañaga}}\email{jcarlos@cimat.mx}

\author*[2]{\fnm{José Luis} \sur{León-Medina}}\email{luis.leon@cimat.mx}

\author[3]{\fnm{Fabiola} \sur{Manjarrez-Gutiérrez}}\email{fabiola.manjarrez@im.unam.mx}

\affil[1]{\orgdiv{Facultad de Ciencias}, \orgname{UNAM},
\orgaddress{\street{Av. Universidad 3000}, \city{Ciudad Universitaria}, \postcode{04510}, \state{CDMX}, \country{México}}}

\affil[2]{\orgdiv{Centro de Investigación en Matemáticas}, \orgname{CIMAT-Mérida}, \orgaddress{\street{Carretera Sierra Papacal, Chuburna Puerto Km 5}, \city{Mérida}, \postcode{97302}, \state{Yucatán}, \country{México}}}

\affil[3]{\orgdiv{Instituto de Matemáticas, Unidad Cuernvaca}, \orgname{UNAM}, \orgaddress{\street{Av. Universidad s/n}, \city{Cuernavaca}, \postcode{62210}, \state{Morelos}, \country{México}}}


\abstract{We propose a way to derive polynomial invariants of closed, orientable $3$-manifolds from Heegaard diagrams via cellularly embedded graphs. Given a Heegaard diagram of an irreducible $3$-manifold $M$, we associate a Heegaard graph $G\subset S$ on the Heegaard surface and restrict to those arising from minimal–genus splittings with a minimal number of vertices. We prove that, up to the natural equivalence of embedded graphs, only finitely many of such minimal Heegaard graphs occur for a fixed manifold. This finiteness enables the definition of $3$-manifold polynomials by evaluating embedded–graph polynomials on representatives of these classes.

For lens spaces we show that the associated Heegaard graphs can be fully classified, and that this classification coincides with the classical one for $L(p,q)$. In this setting the Tutte, Penrose, and Bollobás--Riordan polynomials behave as invariants of lens spaces, and computational evidence suggests that they may in fact be complete invariants. \rev{For the Poincaré homology sphere, whose genus-two Heegaard splitting is unique, we compute the Penrose polynomial of a Heegaard graph with the minimal number of vertices, which we conjecture to be unique; this illustrates that the invariants are computable beyond lens spaces and raises natural questions about broader families of $3$-manifolds.}}

\maketitle

\section{Introduction}\label{sec:intro}

The understanding of $3$-manifolds is a classical topic in low-dimensional topology. Many invariants to distinguish $3$-manifolds have been proposed, ranging from the fundamental group to some sophisticated invariants such as the Heegaard-Floer homology \cite{OS}. One of the most useful representations of $3$-manifolds is the Heegaard diagram, which consists of two families of closed curves on a closed orientable surface and encodes a splitting of the 3-manifold into two handlebodies of the same genus (see \cite{Hempel1976} for an extended exposition). From a Heegaard diagram of a $3$-manifold, one can compute the fundamental group and the Heegaard–Floer homology, among other useful invariants.

In the present work, we propose an approach to compute polynomial invariants of $3$-manifolds. Given a Heegaard diagram of an irreducible $3$-manifold, $M$, we pass to a cellularly embedded Heegaard graph $G$ on the Heegaard surface and focus on those graphs arising from minimal-genus splittings with a minimal number of vertices. Up to the equivalence described in Section~\ref{Preliminaries}, only finitely many of such graphs occur (Theorem~\ref{finiteness}). We then evaluate embedded-graph polynomials on $G$ ---in particular the Tutte, Penrose, and Bollobás–Riordan ribbon-graph polynomials--- to obtain $3$-manifold invariants; computational implementations are available in \cite{FG}. Each one of these polynomials has particular properties, and they have been widely studied \cite{EMM}. This work aims to interpret these polynomial invariants in terms of the topological properties of the underlying $3$-manifolds.

Section~\ref{Preliminaries} recalls Heegaard splittings and diagrams, the passage to cellularly embedded Heegaard graphs, and the equivalence notion between them. Section~\ref{sec:polynomials} reviews the embedded graph polynomials of interest---including the medial-graph construction and the $k$-valuation interpretation---while Section~\ref{sec:3manifoldpolynomials} formulates the resulting $3$-manifold polynomial invariants from representatives of the finite classes of Heegaard graphs.

In the particular case of lens spaces, we show that their Heegaard graphs are circulant graphs of the form $C_p(\pm 1,\pm q)$ and admit a unique embedding class into the Heegaard torus. As a consequence, the classification of lens spaces coincides precisely with the classification of their Heegaard graphs. This structural identification ensures that the resulting $3$-manifold polynomials are computable, and that the associated Tutte, Penrose, and Bollobás–Riordan polynomials are genuine invariants of $L(p,q)$. Moreover, computational evidence strongly suggests that these invariants are complete: while the Tutte and Bollobás–Riordan polynomials already recover the parameter $p$, evaluations also appear to distinguish the orbits of $q$, and the Penrose polynomial provides further structural information, including a characterization of the $q=1$ case (Theorem~\ref{teo:uno}).

Finally, in Subsection~\ref{subsec:Poincare}, we turn to the Poincaré homology sphere as the next natural test case beyond lens spaces. A key step was to determine the minimal number of vertices of its Heegaard graph, which turns out to be $12$. This observation, originally communicated to us by Francisco González-Acuña (see Appendix), made it possible to carry out explicit computations of the proposed polynomials. 

\section{Preliminaries}\label{Preliminaries}

\subsection{Heegaard diagrams}\label{subsec:Heegaard}

Throughout the paper let $M$ denote a connected, closed, orientable
$3$--manifold. 

\begin{figure}
    \centering
    \includegraphics[trim=150pt 172pt 50pt 110pt, clip, width=0.68\linewidth]{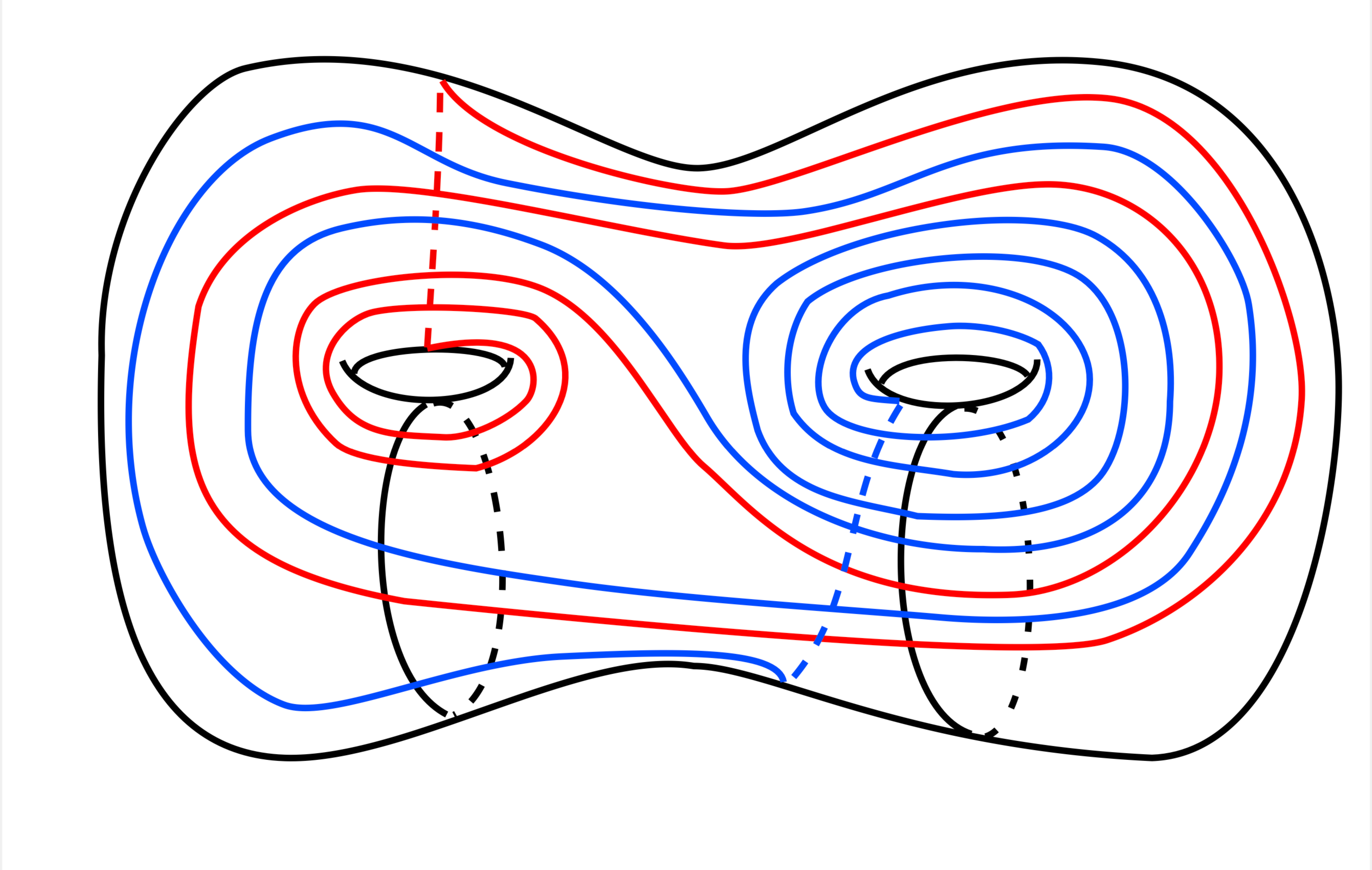}
    \caption{A genus–two Heegaard diagram for the Poincaré sphere $P_3$ with 12 vertices}
    \label{fig:PoincareSphere}
\end{figure}

\begin{definition}
A \emph{Heegaard splitting} of $M$ is a decomposition
$$
M \;=\; V \cup_S W ,
$$
where $V$ and $W$ are handlebodies of the same genus $g\ge 0$ and
$S = \partial V = \partial W$ is their common boundary, the
\emph{Heegaard surface}.  
The integer $g$ is called the \emph{genus of the splitting}.

From this splitting, we can choose pairwise disjoint, properly embedded \emph{meridian disks}
$$
D_V = \{D^V_1,\dots,D^V_g\}\subset V, 
\qquad
D_W = \{D^W_1,\dots,D^W_g\}\subset W,
$$
such that cutting along them turns each handlebody into a $3$–ball. Let
$$
v = \{\partial D^V_1,\dots,\partial D^V_g\},
\qquad
w = \{\partial D^W_1,\dots,\partial D^W_g\}
$$
be the corresponding collections of essential simple closed curves on
$S$.  
The triple
$
(S; v, w),
$
so obtained is a Heegaard diagram associated with the splitting. 

More generally, a \emph{Heegaard diagram} can be defined independently as a triple
$
\mathcal D = (S; v, w),
$
where $S$ is a closed orientable surface of genus $g$, while $v$ and $w$ are two collections, each consisting of $g$ pairwise disjoint simple closed curves, such that the surfaces $S-v$ and $S-w$ are connected.
Any such diagram provides enough information to recover a Heegaard splitting, and therefore, a $3$-manifold \cite[Remark 6.2.2.]{Schultens2014p}.

We will always assume that in the Heegaard diagrams, the curves $v_i\in v$ and $w_j\in w$ intersect each other transversally. For later arguments, it is convenient to identify $S$ with a
\emph{standard genus-$g$ surface} and regard the curves of $v$ as the {standard} meridian system of $S$.
Under this convention, a Heegaard diagram is completely specified by the
second set of curves $w=\{w_1,\dots ,w_g\}$; see
Fig.~\ref{fig:PoincareSphere} for the case of the Poincaré homology sphere.
\end{definition}

\rev{\begin{definition}\label{def:diagram-equiv}
Two Heegaard diagrams of a $3$-manifold $M$ are \emph{equivalent} if their underlying Heegaard splittings are equivalent, that is, if some self-homeomorphism of $M$ carries one splitting onto the other.
\end{definition}}

\rev{\begin{remark}
It is known that this is equivalent to passing from one diagram to the other by homeomorphisms, semi-isotopies, and sums of meridians; see \cite[Ch.~5]{Fomenko1997p}.
\end{remark}}

\subsection{Cellularly embedded graphs}

Heegaard splittings and their associated diagrams provide a useful way of decomposing $3$-manifolds. However, a significant disadvantage is that the splitting for a given $3$-manifold $M$ is not unique; in fact, there are infinitely many such splittings for $M$. To motivate our next steps, we collect the basic notions of finite graphs, their embeddings in surfaces, and the ribbon graphs that arise
from such embeddings; these concepts will let us speak conveniently about the graph polynomials associated to a cellularly
embedded graph. We adopt the terminology of
\cite{GrossTucker87}.

\begin{definition}
    A \emph{finite graph} $G$ is a pair $(V, E)$, where $V$ is a finite set of \emph{vertices} and $E$ is a set of \emph{edges}. Each edge $e \in E$ is associated with a set of one or two vertices from $V$, its \emph{endpoints}. We write $V(G):=V$ and $E(G):=E$ to refer to the vertex and edge sets of $G$. \rev{The \emph{geometric realization} of $G$, denoted $|G|$, is the $1$-dimensional CW complex with one $0$-cell for each vertex and one $1$-cell for each edge, the latter attached to the $0$-cells corresponding to its endpoints.}
    
    A \emph{walk} is a sequence $v_0, e_1, v_1, \dots, e_n, v_n$ where for each $j$, the vertices $v_{j-1}$ and $v_j$ are the endpoints of the edge $e_j$. The walk is \emph{closed} if $v_0 = v_n$ and is a \emph{path} if all its vertices are distinct. A graph is \emph{connected} if any pair of its vertices can be joined by a path.
    
    We denote the number of vertices, edges, and connected components of $G$ by $v(G)$, $e(G)$, and $k(G)$, respectively. The \emph{rank} of $G$ is $r(G) := v(G) - k(G)$, and its \emph{nullity} is $n(G) := e(G) - r(G)$. {A subgraph $H=(U,A)$ of $G=(V,E)$ is a \emph{spanning subgraph} if $U=V$. If $U$ consists of exactly the vertices of $G$ that are incident to edges in $A\subset E$, we will say that $H$ is the \emph{subgraph induced} by $A$ and denote it by $G|_A$. Given a graph $G=(V,E)$ and $A\subset E$, we let $e(A)$, $v(A)$, $k(A)$, $r(A)$, and $n(A)$ denote the number of edges, vertices, components, rank, and nullity of the spanning subgraph $(V,A)$ of $G$.} \rev{These quantities are used later in Definition~\ref{def:BR}, where the Bollobás--Riordan polynomial is defined through the exponents $r(G)-r(A)$ and $n(A)$.}

    {A \emph{graph map} $\phi:G\to G'$ consists of a vertex function $\phi_V: V(G) \to V({G'})$ and an edge function $\phi_E: E(G)\to E(G')$ such that for each $e \in E(G)$, $V(\phi_E(e))= \phi_V(V(e))$. A graph map is an \emph{isomorphism} if both of its vertex function and its edge function are bijective. Two graphs $G_1$ and $G_2$ are \emph{isomorphic}, denoted $G_1 \cong G_2$, if there is an isomorphism from one to the other.}
\end{definition}

\begin{definition} Let $F$ be a compact connected surface without boundary. An \emph{embedding} of a graph $G$ in $F$ is a continuous injective map $\iota$ from the geometric realization of $G$, $|G|$, into the surface.
The embedded graph $\iota(G)$ is said to be \emph{cellular} if every
connected component of $F \setminus \iota(G)$ is homeomorphic to an
open disk. These components are called the \emph{faces} of the embedding, and their number is denoted by $f(G)$.

By abuse of terminology we identify a graph $G$ with its image in $F$---the embedded graph---and use ‘abstract graph’ to refer to the original
combinatorial object $G$. 

We will say that two cellularly embedded graphs $G\subset F$ and $G'\subset F'$ are \emph{equivalent} if there is a homeomorphism $\varphi: F\to F'$ such that $\varphi(G)=G'$ and $\varphi\big\lvert_G:G\to G'$ is an isomorphism of abstract graphs.
\end{definition}

\rev{The definition above does not require $F$ to be orientable. Throughout this paper, however, the surface $F$ is always orientable, as it will be a Heegaard surface; we keep the general formulation only to fix the terminology.}

\rev{Cellular embeddings in orientable surfaces are encoded combinatorially by rotation systems.}

\rev{\begin{definition}\label{def:rotation}
Each edge of a graph $G$ consists of two \emph{darts} (or half-edges); let $D$ be the set of all darts of $G$, so that $|D|=2\,e(G)$. A \emph{rotation system} (or \emph{combinatorial map}; see~\cite{LandoZvonkin2004}) on $G$ is a pair $(\sigma,\rho)$ of permutations of $D$, where $\rho$ is a fixed-point-free involution whose orbits are the edges of $G$, and the orbits of $\sigma$ are the vertices, each recording the cyclic order of the darts around that vertex. The map is \emph{$4$-regular} when every orbit of $\sigma$ has four darts. Replacing $\sigma$ by $\sigma^{-1}$ reverses the orientation and yields the \emph{symmetric} rotation system $(\sigma^{-1},\rho)$.
\end{definition}}

\rev{\begin{theorem}[Heffter--Edmonds; see \cite{MoharThomassen2001}]\label{prop:ribbon-rotation}
For a connected graph $G$, the cellular embeddings of $G$ into closed orientable surfaces, taken up to orientation-preserving homeomorphism, are in bijection with the rotation systems on $G$; under this bijection the faces of the embedding are the cycles of $\varphi=\rho\sigma$.
\end{theorem}}

\rev{This combinatorial encoding is the one we use to compare the cellularly embedded graphs arising from Heegaard diagrams and to prove finiteness (Theorem~\ref{finiteness}). Equivalently, a cellular embedding can be described by thickening the graph within its surface to a ribbon graph, the natural setting for the Bollobás--Riordan and Penrose polynomials.}

\rev{\begin{definition}\label{def:ribbon}
A \emph{ribbon graph} $G=(V(G),E(G))$ is a (possibly non-orientable) surface with boundary, represented as the union of a set $V(G)$ of \emph{vertex-disks} and a set $E(G)$ of \emph{edge-bands}, where each edge-band is a topological rectangle, such that:
\begin{enumerate}[label=\textup{(\arabic*)}, nosep]
\item the vertex-disks and edge-bands intersect in disjoint line segments;
\item each such line segment lies on the boundary of precisely one vertex-disk and precisely one edge-band; and
\item every edge-band contains exactly two such line segments.
\end{enumerate}
In this sense, the vertex-disks may be viewed as $0$-handles and the edge-bands as $1$-handles attached along two disjoint boundary intervals.
The ribbon graph is \emph{orientable} if the underlying surface is orientable.
\end{definition}}

\subsection{Embedded graph polynomials}\label{sec:polynomials}

Once the concept of cellularly embedded graph has been defined we will discuss several generalizations of classical graph polynomials to this new context. We are particularly interested in the Tutte, Penrose and Bollobás-Riordan polynomials. For an extended exposition on these polynomials, refer to \cite{EMM}. 

The Tutte polynomial \cite{T} is a classical $2$-variable polynomial in graph theory that generalizes other polynomials such as the chromatic polynomial or the flow-polynomial for abstract graphs. This polynomial has information on the connectedness of the graph and is related to some physical models such as the Potts model. A generalization of the Tutte polynomial for embedded graphs is the Bollobás-Riordan ribbon graph polynomial, introduced in \cite{Bo}. The Bollobás-Riordan polynomial has four variables and contains topological information on the subgraphs of a ribbon graph, such as the number of connected components and the number of boundary components. The precise formulation is as follows.

\begin{definition}\label{def:BR}
Let $G$ be a cellularly embedded graph. The \emph{Bollobás-Riordan polynomial} $R(G;x,y,z,w)$ is defined as
\[
R(G;x,y,z,w) = \sum_{A\subset E(G)} (x-1)^{r(G)-r(A)}y^{n(A)}z^{k(A)-f(A)+n(A)}w^{t(A)}.
\]
where the sum runs over all spanning ribbon subgraphs $G|_A$ obtained from every possible edge subset $A\subset E(G)$; $f(A)$ denotes the number of boundary components (also called faces) 
of $G|_A$ in the associated ribbon graph; and $t(A)=0$ if $G|_A$ is orientable and $t(A)=1$ otherwise.
\end{definition}

\begin{definition}
Let $G$ be a cellularly embedded graph. The \emph{Tutte polynomial} $T(G;x,y)$ is defined as $T(G;x,y)=R(G;x,y-1,1,1).$
\end{definition}
    
Similarly, the Penrose polynomial was introduced for planar graphs by Penrose \cite{P} in 1971 and extended in 2013 by  Ellis-Monaghan and Moffat to embedded graphs \cite{ELLISMONAGHAN2013424}. It encodes many combinatorial properties of the embedded graph, including the number of some colourings of the graph, as we will see below.

\begin{figure}[htb]
  \centering
    \includegraphics{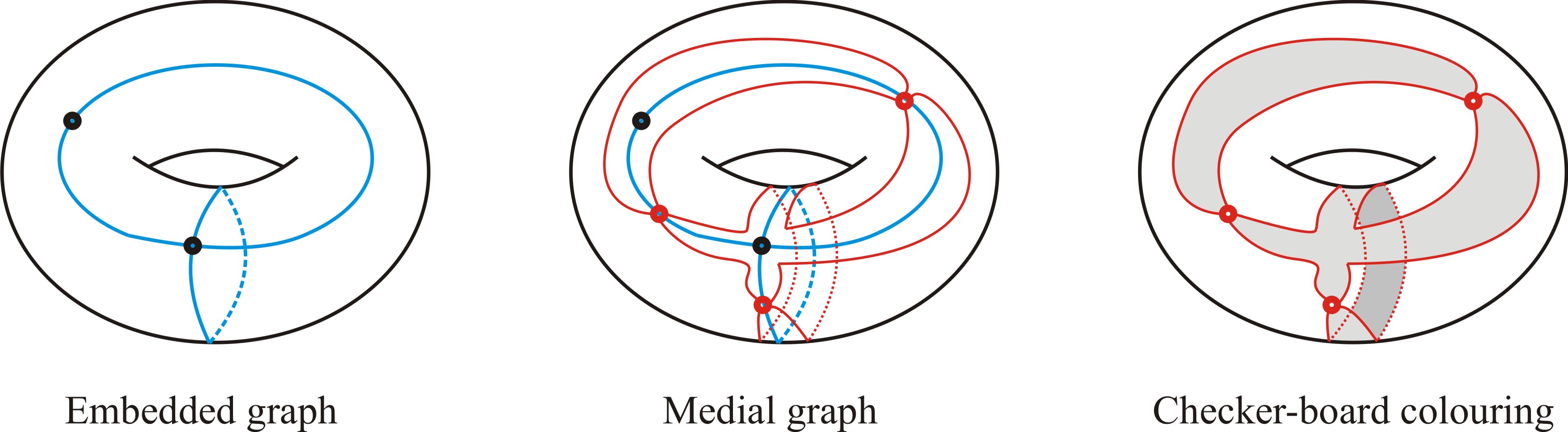}
  \caption{Construction of the medial graph}
  \label{fig:f4}
\end{figure}

Given a cellularly embedded graph $G$, its \emph{medial graph} $G_m$ is
constructed by placing a vertex of degree $4$ on each edge of $G$ and
joining these new vertices by arcs that follow the face boundaries of
$G$, see Fig. ~\ref{fig:f4}. The medial graph $G_m$ is itself
cellularly embedded in the same surface as $G$, and it is always
$4$--regular. Moreover, $G_m$ admits a canonical \emph{checkerboard
colouring} of its faces: the faces incident with the vertices of $G$ are
coloured black, and the remaining faces are coloured white as shown on the right side of Fig. ~\ref{fig:f4}.

At each vertex of $G_m$ there are three possible \emph{vertex states},
obtained by merging in pairs the four incident edges with the vertex. Using the
checkerboard colouring, these states are distinguished as:
\begin{itemize}
\item the \emph{white smoothing}, where the paired edges connect through
the white faces,
\item the \emph{black smoothing}, where they connect through the black
faces, and
\item the \emph{crossing}, where opposite edges are connected in a crossing fashion.
\end{itemize}

\rev{Here a crossing is a purely local recombination of the four half-edges meeting at a vertex of $G_m$: it is neither an additional vertex nor a knot-theoretic crossing (there is no over/under information). It simply joins the two pairs of opposite edges transversally, which amounts to attaching the corresponding band with a half-twist. Consequently the ribbon graph determined by a state is an abstract ribbon graph in the sense of Definition~\ref{def:ribbon} and need not be embedded in $F$; this causes no difficulty, since only its number of boundary components enters the Penrose polynomial.}

Given a cellularly embedded graph $G$ and its medial graph $G_m$, a \emph{Penrose state} of $G_m$ is a choice of white smoothing or crossing at every vertex of $G_m$. Note that the black smoothing, although available in the general notion of vertex states, is not used in the definition of the Penrose polynomial. Each such state $s$ determines a ribbon graph obtained from the black faces of $G_m$, whose vertices coincide with the vertices of $G$, while the ribbon edges may be twisted with respect to the ribbon graph associated to $G$.

For a Penrose state $s$, we write $c(s)$ for the number of boundary components of the ribbon graph associated to $s$, and $cr(s)$ for the number of vertices of $G_m$ at which $s$ takes the crossing state.
    
\begin{definition}
\label{df:pen} 
Let $G$ be a cellularly embedded graph in a closed surface $F$, and let
$G_{m}$ be its medial graph. The \emph{Penrose polynomial} of $G$ is
\[
P(G;\lambda) \;=\; \sum_{s\in \mathcal{P}(G_{m})} (-1)^{\,cr(s)}\,\lambda^{\,c(s)} ,
\]
where $\mathcal{P}(G_{m})$ denotes the set of Penrose states of $G_m$.
\end{definition}

\begin{figure}[hbt]
  \centering
    \includegraphics{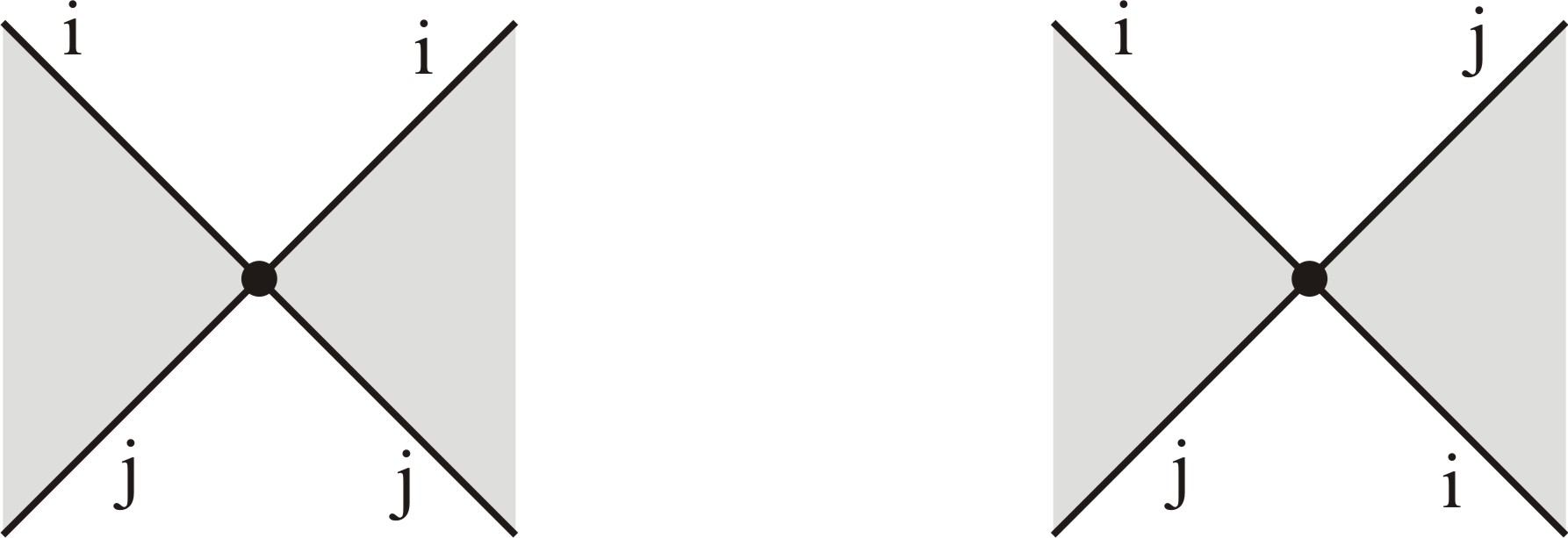}
  \caption{Admissible $k$-valuation}
  \label{fig:f6}
\end{figure}

The Penrose polynomial of a graph $G$ is related to the $k$-valuations of the medial graph $G_m$. Suppose that $G_m$ is canonically checker-board coloured. For a natural number $k$, a \emph{$k$-valuation} of $G_m$ is an edge $k$-colouring $\phi: E(G_m)\rightarrow \{1,2,\ldots, k\}$ such that for each $1\leq i \leq k$ and each vertex of $G_m$, the number of $i$-coloured edges of $G_m$ incident with the vertex is even.  We say that a $k$-valuation is \emph{admissible} if at each vertex of $G_m$ the $k$-valuation is one of the two types in Fig. ~\ref{fig:f6}, where $i\neq j$. Note that these two types of colouring at a given vertex correspond to the white smoothing and the crossing state of the vertex. We next recall a theorem connecting the integer evaluations of the Penrose polynomial of a graph with the admissible $k$-valuations of the corresponding medial graph.    

\begin{theorem}[\protect{\cite[Thm. 6.8]{ELLISMONAGHAN2013424}}]\label{teo:val}
Let $G$ be an embedded graph and take $k\in \mathbb{N}$, then $P(G;k)=\sum (-1)^{cr(s)}$, where the sum is taken over all admissible $k$-valuations $s$ of $G_m$ and $cr(s)$ is the number of crossing states in $s$.  
\end{theorem}

\section{Heegaard graphs}\label{heegaardgraphs}

In what follows we will establish a connection between $3$-manifolds and cellularly embedded graph polynomials via the determination of a finite number of cellularly embedded graphs associated to the Heegaard splittings of a given $3$-manifold. This section is devoted to the construction of such graphs.

\begin{definition}
Let $M$ be a 3-manifold and $V\cup_S W$ a Heegaard splitting for $M$, of genus $g$. Let  $\mathcal D=(S; v,w)$ be a Heegaard diagram for  $V\cup_S W$. The  \emph{Heegaard graph} associated with $\mathcal D$ is the embedded
$1$–complex
$$
G_{\mathcal D}= v\cup w
$$
The vertices and edges of the abstract graph are defined as follows:
\begin{itemize}
\item The vertex set is the set of intersection points
      $$
      V\bigl(G_{\mathcal D}\bigr)= v\cap w .
      $$
\item There is an edge joining $p$ and $q$ for each connected component of $G_{\mathcal D}\setminus V(G_{\mathcal D})$ having $p$ and $q$ as boundary points.
\end{itemize}

When no confusion is likely, we identify $G_{\mathcal D}$ with its image {in} $S$ and call it simply the \emph{Heegaard graph} of $\mathcal{D}$.
\end{definition}

{The graph associated with an arbitrary Heegaard diagram is not guaranteed to be cellularly embedded. Furthermore, a  fundamental problem arises if a curve in one system is disjoint from the other, as the combinatorial structure of the graph degenerates. We will show that both of these issues can be avoided by imposing conditions on the manifold $M$ and on the minimality of $|v\cap w|$.}

\begin{definition}
    A $3$-manifold $M$ is \emph{irreducible} if each $2$-sphere in $M$ bounds a $3$-ball in $M$.
\end{definition}

\begin{definition}
    Let $M=V\cup _S W$ be a Heegaard splitting and let $\mathbb{S}^3= V'\cup_T W'$ be the standard genus 1 Heegaard splitting of $\mathbb{S}^3$. The pairwise connected sum $(M,S)\sharp (\mathbb{S}^3, T)$ defines a Heegaard splitting $M=\tilde{V}\cup_{\tilde{S}}\tilde{W}$ called an \emph{elementary stabilization} of $M=V\cup_S W$. A Heegaard splitting is called a \emph{stabilization} of $M=V\cup_SW$ if it is obtained from $M=V\cup_SW$ by performing a finite number of elementary stabilizations.
\end{definition}

\begin{definition}
    A Heegaard splitting $M=V\cup _S W$ is \emph{reducible} if there is a 2-sphere $\Sigma$ in $M$ such that $\Sigma \cap S$ is an essential closed curve. A Heegaard splitting is \emph{irreducible} if it is not reducible. 
\end{definition}

\begin{theorem}[\cite{Scharlemann}]\label{theo:redicibleimpliesstab}
    Suppose $M$ is an irreducible 3-manifold and $V\cup _S W$ is a reducible Heegaard splitting of $M$. Then $V\cup _S W$ is stabilized.
\end{theorem}

From now on, we will consider irreducible manifolds with irreducible Heegaard splittings

\begin{proposition}\label{prop:irredHeegSplitiscellular}
    Let $V\cup_S W$ be an irreducible Heegaard splitting, and let $(S; v, w)$ be any Heegaard diagram of $V\cup_S W$. Then $S\setminus \{v,w\}$ is a collection of open disks.
\begin{proof}
    Suppose $S\setminus \{v,w\}$ contains a component $C$ which is not a disk. Then there is an essential curve $\beta$ in  $C$ that is disjoint from $v\cup w$. Cutting $V$ along $v$ gives a 3-ball, since $\beta$ is disjoint from $v$, then $\beta$ bounds a disk in such a 3-ball. After reglueing to recover $V$, the curve $\beta$ will bound a disk in $V$. By analogous argument, $\beta$ bounds a disk \rev{in} $W$. These two disks define a 2-sphere which intersects $S$ in an essential curve $\beta$, which contradicts the irreducibility of $V\cup_S W$. Therefore all components of $S\setminus \{v,w\}$ are open disks. 
\end{proof}
\end{proposition}

To work with a more tractable class of diagrams, we impose two successive minimality conditions. It is known that any connected, closed, orientable $3$-manifold has a Heegaard splitting. First, we choose a Heegaard splitting of minimal genus.  

\begin{definition}
Let $M$ be any closed, connected, orientable $3$-manifold. The \emph{Heegaard genus} of $M$, denoted $h(M)$, is the minimum genus of any Heegaard splitting for~$M$. A Heegaard splitting realizing $h(M)$ is called a \emph{minimal Heegaard splitting}
\end{definition}

Next, among such minimal Heegaard splittings, we select those whose associated graphs have the fewest possible vertices.

\begin{definition}
Let $\mathcal{H}(M)$ be the set of all Heegaard graphs associated with diagrams of genus $h(M)$. Let $\mathcal{G}(M) \subset \mathcal{H}(M)$ be the sub-collection of those graphs having the minimal possible number of vertices.
\end{definition}

Our next step toward defining a potential polynomial invariant is to show that the set $\mathcal{G}(M)$ consists of only a finite number of graph types. More formally:

\begin{theorem}\label{finiteness} Let $M$ be irreducible. The set of equivalence classes of $\mathcal{G}(M)$ is finite.
\begin{proof}
\rev{If $h(M)=0$ then $M=S^3$ and $\mathcal{G}(M)$ consists of a single, empty graph, so the statement is trivial; we therefore assume $h(M)>0$, in which case every minimal Heegaard graph has $m\geq 1$ vertices and is $4$-valent.}

A minimal Heegaard splitting is irreducible. Hence by Proposition \ref{prop:irredHeegSplitiscellular} any Heegaard graph, $G_\mathcal{D}$ with a minimal number of vertices, say $m$, is cellularly embedded in the splitting surface $S_g$.
\rev{Being $4$-valent, $G_\mathcal{D}$ has $e=2m$ edges and, by Euler's formula, exactly $d=\chi(S_g)+m$ faces. By Theorem~\ref{prop:ribbon-rotation}, the cellular embedding $G_\mathcal{D}\subset S_g$ is determined, up to homeomorphism, by a rotation system on its abstract graph. Viewed in this way, the equivalence classes of $\mathcal{G}(M)$ form a subset of the $4$-regular rotation systems on $m$ vertices with $d$ faces, taken up to isomorphism and up to the symmetric rotation system $\sigma\mapsto\sigma^{-1}$. This set is finite: there are finitely many $4$-valent graphs on $m$ vertices, and on each the involution $\rho$ is fixed while $\sigma$ is a choice of cyclic order of the four darts at every vertex, giving at most $(3!)^{m}=6^{m}$ rotation systems. Therefore $\mathcal{G}(M)$ has only finitely many equivalence classes.}
\end{proof}
\end{theorem}

\rev{We stress that Theorem~\ref{finiteness} is purely an existence (finiteness)
statement: it does not provide an effective procedure to determine the elements
of $\mathcal{G}(M)$. Even computing the Heegaard genus is NP-hard, and the
corresponding algorithms are primarily existence results, with no general
implementations currently available \cite{BurtonThompson2025}. The effective
determination of $\mathcal{G}(M)$ thus lies outside the scope of the present
work.}

\rev{We further note that, although only finitely many $4$-valent graphs with $m$
vertices exist, their number grows rapidly with $m$. Characterizing precisely
which $4$-valent graphs are realized by a Heegaard diagram of some $3$-manifold
is an interesting question, which we leave open: not every $4$-valent graph
arises in this way.}

\section{3-manifold polynomials}\label{sec:3manifoldpolynomials}

Tutte polynomial has been implemented in many known libraries such as \emph{Sage} in \emph{Python}. In \cite{FG}, the authors implemented computer programs to compute the Penrose and Bollobás--Riordan polynomials based on the state formulations for these polynomials. These implementations were very useful to gain insight into the nature of the polynomials of Heegaard graphs of lens spaces.

\subsection{The case of lens spaces}\label{lens}

Lens spaces are genus-$1$ Heegaard splittings. That is to say, a Heegaard diagram of a lens space is represented by a standard meridian and a simple closed curve on a torus, or a torus knot; and this curve can be parametrized by two relative prime numbers $p$ and $q$. Lens spaces are fully classified in terms of these two numbers according to the following conditions (see \cite{Br}):      

\begin{theorem}[Classification of lens spaces]\label{teo:clas}
Two lens spaces $L(p,q)$ and $L(p',q')$ are homeomorphic if and only if
\begin{itemize}\itemindent=1.5em
\item[(i)] $p=p'$, and
\item[(ii)] $q'\equiv \pm q^{\pm 1} \pmod p$.
\end{itemize}
\end{theorem}

It turns out that the kind of graphs obtained as Heegaard graphs of lens spaces are a special type of a well-known family of circulant graphs. We recall the following definition.

\begin{definition}
Let $n\ge 1$ and let $S\subseteq \mathbb{Z}_n\setminus\{0\}$.  
The \emph{circulant} graph $C_n(S)$ has vertex set $\mathbb{Z}_n$, and two vertices  
$v_i,v_j$ are adjacent if and only if $i-j\pmod n\in S$.
\end{definition}

\begin{remark}\label{rem:lens-circulant}
For the standard genus-$1$ Heegaard diagram of $L(p,q)$, the associated Heegaard graph is the circulant $C_p(\pm 1,\pm q)$: its vertices are the $p$ intersection points around the torus meridian, and edges come from the two families of arcs joining $i$ with $i\pm 1$ and $i\pm q$ (mod $p$).
\end{remark}

While the embedded Heegaard graph carries strictly more information than its abstract counterpart, it turns out that for lens spaces this additional structure is unnecessary. The abstract isomorphism type of the Heegaard graph alone suffices to recover the homeomorphism type of the lens space, showing that Heegaard graphs provide a complete invariant in this setting. The following theorem establishes this fact.

\begin{theorem}\label{thm:lens-circulant-abstract}
For integers $p\ge 3$ and $q,q'$ coprime with $p$, 
the lens spaces $L(p,q)$ and $L(p,q')$ are homeomorphic 
if and only if their Heegaard graphs $C_p(\pm 1,\pm q)$ and $C_p(\pm 1,\pm q')$ are isomorphic as abstract graphs. 

\begin{proof}
It is known that the isomorphism problem for 
double--loop circulant graphs, i.e.\ circulants with connection set 
$S=\{\pm a,\pm b\}$ of size four, has a complete solution: 
every such set $S$ is a CI--subset of $\mathbb{Z}_p$, so that
\[
C_p(S)\cong C_p(S') \iff S' = aS \quad \text{for some } a\in \mathbb{Z}_p^\times.
\]
See \cite[Thm.~5.4]{muzychuk}.

Specializing to $S=\{\pm 1,\pm q\}$ and $S'=\{\pm 1,\pm q'\}$, 
this means
\[
C_p(\pm 1,\pm q)\cong C_p(\pm 1,\pm q') \iff 
\{\pm 1,\pm q\} = \{\pm a,\pm aq'\}
\quad \text{for some } a\in \mathbb{Z}_p^\times.
\]

We now split into two cases:
\begin{itemize}
\item If $a\equiv \pm 1 \pmod p$, then $q'\equiv \pm q \pmod p$.
\item If $a\equiv \pm q \pmod p$, then 
$\{\pm 1,\pm q\}=\{\pm q,\pm(qq')\}$, which forces 
$qq'\equiv \pm 1\pmod p$, hence $q'\equiv \pm q^{-1}\pmod p$.
\end{itemize}

Thus $C_p(\pm 1,\pm q)\cong C_p(\pm 1,\pm q')$ implies 
$q'\equiv \pm q^{\pm 1}\pmod p$, exactly as in 
Theorem~\ref{teo:clas}, so $L(p,q)\cong L(p,q')$. 
Conversely, if $q'\equiv \pm q^{\pm 1}\pmod p$, then 
choosing $a\in\{\pm 1,\pm q'\}$ with $S'=aS$ yields 
$C_p(\pm 1,\pm q)\cong C_p(\pm 1,\pm q')$. 
\end{proof}
\end{theorem}

Regarding the construction of Heegaard graphs from Section~\ref{heegaardgraphs}, in the case of lens spaces, since there exists essentially a unique splitting disk in each of the two solid tori, we have:

\begin{remark}
Given a lens space $L(p,q)$, there exists a unique embedded graph obtained from the standard Heegaard diagram of  $L(p,q)$. 
\end{remark}

Moreover, the Tutte, Penrose and ribbon graph polynomials computed for this graph are polynomial invariants of the corresponding lens space as we will show.

\begin{figure}
  \centering
    \includegraphics{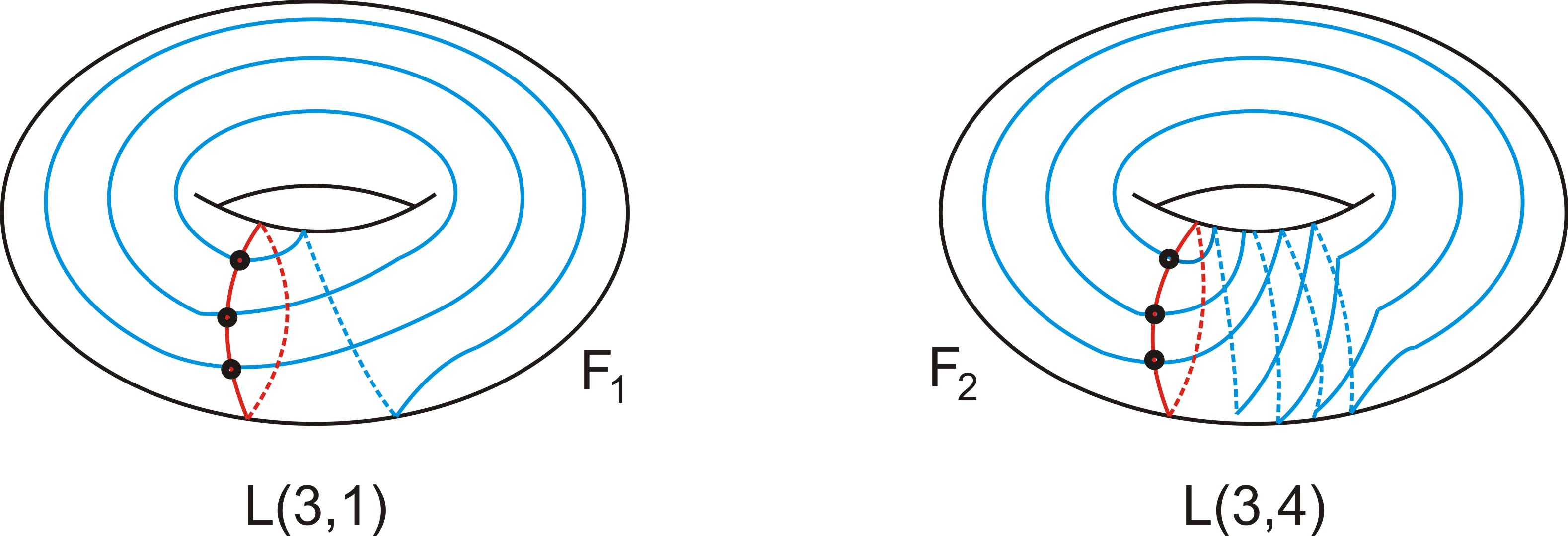}
  \caption{Homeomorphic lens spaces}
  \label{fig:f1}
\end{figure}

\begin{figure}
  \centering
    \includegraphics{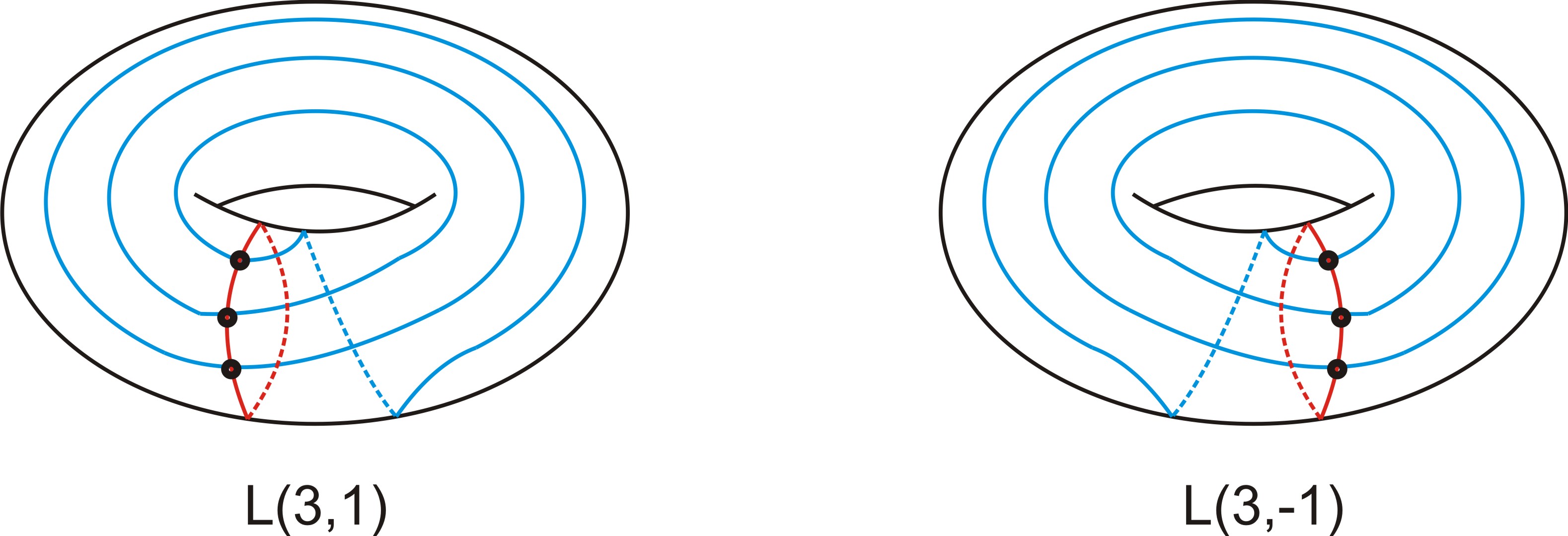}
  \caption{Homeomorphic lens spaces}
  \label{fig:f2}
\end{figure}

\begin{theorem}\label{teo:inv}
If $L(p,q)$ and $L(p,q')$ are homeomorphic lens spaces with embedded Heegaard graphs $G_1$ and $G_2$, respectively, then $G_1\subset F_1$ and $G_2\subset F_2$ are equivalent cellularly embedded graphs with respect to the Heegaard tori $F_1$ and $F_2$.

\begin{proof}
According to the lens spaces classification,  $q'\equiv \pm q^{\pm 1} \mod p$. We analyze the following cases.
\begin{itemize}
\item[]\emph{Case 1:} $q\equiv q' \mod p$.\\
Without loss of generality, we can assume that $0\leq q< p$, then $q'=p\cdot r+q$ for some $r\in \rev{\mathbb{Z}}$. Let $\varphi_{r}:F_1\rightarrow F_2$ be the homeomorphism given by $r$ Dehn twists along a meridian curve in $F_1$ (see Fig. \ref{fig:f1}). The homeomorphism $\varphi_{r}$ induces an isomorphism of the cellularly embedded graphs  $G_{1}$ and $G_2$. According \rev{to} this fact, we can assume that $0\leq |q'| < p$.
\item[]\emph{Case 2:} $q'\equiv -q\mod p$.\\
Consider the homeomorphism $\varphi: F_1\rightarrow F_2$ given by a reflection along a plane (Fig. ~\ref{fig:f2}). This homeomorphism induces an equivalence of the embedded graphs $G_1$ and $G_2$.
\item[]\emph{Case 3:} $q'\equiv q^{-1} \mod p$.\\ 
If $H_1$  and $H_2$ are the two tori in the Heegaard decomposition of $L(p,q)$, there exists a homeomorphism $\psi: L(p,q)\rightarrow L(p,q')$ that swaps the tori $H_1$ and $H_2$. Such homeomorphism induces a homeomorphism between $F_1$ and $F_2$ which induces an equivalence of the embedded graphs $G_1$  and $G_2$.     
\end{itemize}  
\end{proof}
\end{theorem} 

As a consequence of Theorems~\ref{thm:lens-circulant-abstract} and \ref{teo:inv}, we get the following corollary.

\begin{corollary}[Classification via Heegaard graphs]\label{cor:lens-graphs}
Two lens spaces $L(p,q)$ and $L(p,q')$ are homeomorphic if and only if
their Heegaard graphs are equivalent.
\end{corollary}

\subsubsection{Lens spaces polynomials}

The computational experiments in \cite{FG} suggest that the ability of
Heegaard graphs to distinguish lens spaces may already be retained by
certain polynomial invariants. In other words, it seems plausible that
one does not need the full Heegaard graph to recover the
homeomorphism type of a lens space, but that the evaluation of suitable
graph polynomials could already yield complete invariants. In this
subsection we present partial results in this direction, together with
open questions relating lens spaces, graph polynomials, and spectral
graph theory.

\paragraph{The Tutte and Bollobás--Riordan polynomials}

Since the Tutte polynomial is a specialization of the Bollobás--Riordan polynomial we show how the parameters $p$ and $q$ are related to the Tutte polynomial.

\begin{proposition}\label{teo:tutte}
Let $G$ be the embedded Heegaard graph corresponding to $L(p,q)$. In $T(G;x,y)$, the coefficient of the highest power of $y$ is $1$ and the corresponding exponent is $p+1$.
\begin{figure}
  \centering
    \includegraphics{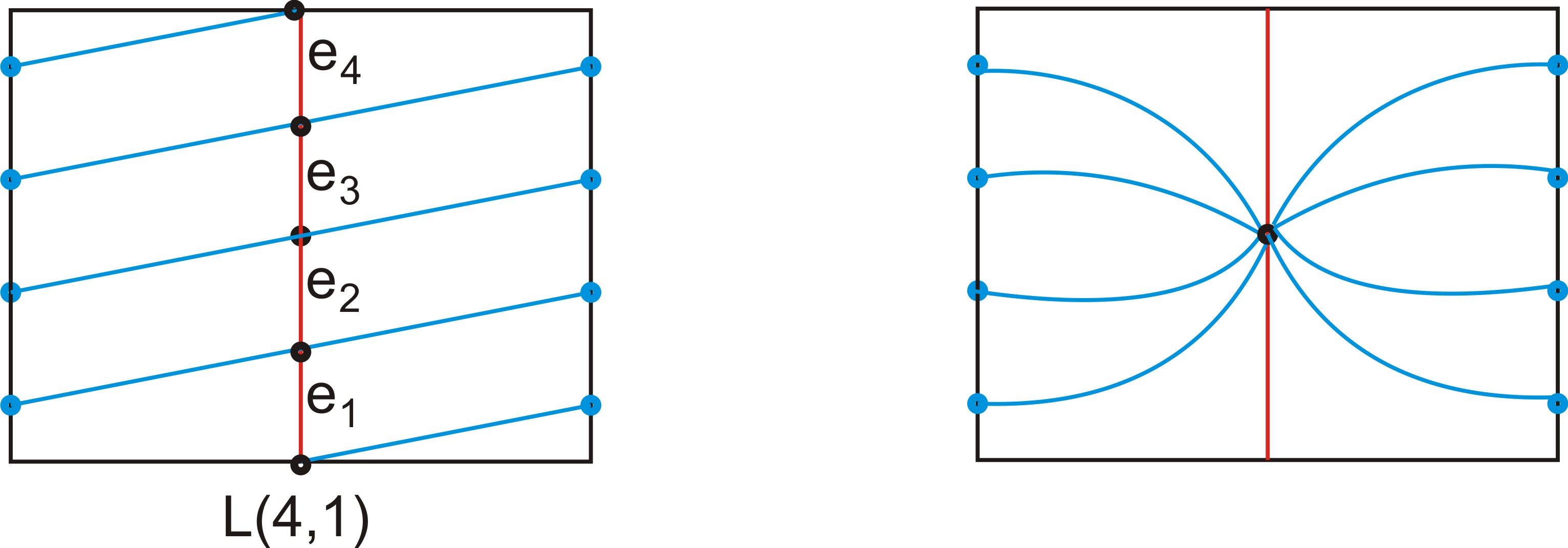}
  \caption{Recursive Tutte polynomial}
  \label{fig:f3}
\end{figure}
\begin{proof}
Order the edges of $G$, $e_1,\dots,e_p$, along the boundary of one meridian disk. Fig.~\ref{fig:f3} shows the Heegaard torus as a square with identified sides. Apply the deletion–contraction recurrence $T(G;x,y)=T(G/e;x,y)+T(G-e;x,y)$ successively on $e_1,\dots,e_p$. The terminal stage with the maximal number of loops (right-hand side diagram of Fig.~\ref{fig:f3}) contributes to the polynomial with the monomial $y^{p+1}$ with coefficient $1$.
\end{proof}
\end{proposition}

Hence, the $p$ value is recoverable from the Tutte and Bollobás--Riordan 
polynomials, but the $q$ parameter exhibits a more subtle behavior. 
Recall that the Tutte polynomial evaluated at $(1,1)$ gives the number 
of spanning trees of the graph  \cite[Prop. 3.1]{ellis-monaghan2022}. This motivates us to focus on the 
spanning tree count of the Heegaard graph $C_p(\pm1,\pm q)$, denoted 
$\tau(p,q)$. This quantity admits a closed product formula via the 
Matrix--Tree Theorem \cite{ZHANG2000337}, and its analysis naturally involves tools from 
spectral graph theory. We will therefore adopt this combinatorial 
viewpoint, referring to standard references such as Biggs’ 
\cite{biggs1974} and Godsil--Royle’s \cite{godsil2001} textbooks, 
as well as the explicit results of Zhang--Yong--Golin \cite{ZHANG2000337}.

\begin{theorem}\label{thm:product-and-parity}
Let $G=C_p(\pm 1,\pm q)$ with $\gcd(p,q)=1$, $p\ge3$. Then
\[
\tau(p,q)=\frac{1}{p}\prod_{j=1}^{p-1}\Bigl(4-2\cos\tfrac{2\pi j}{p}-2\cos\tfrac{2\pi q j}{p}\Bigr),
\]
and $\tau(p,q)$ depends only on the orbit $\{\pm q^{\pm1}\}\subset(\mathbb{Z}/p\mathbb{Z})^\times$.
\end{theorem}

\begin{proof}
By \cite[Lemma~1]{ZHANG2000337}, for an undirected circulant $C_n(s_1,\dots,s_k)$ we have
\[
T\bigl(C_n(s_1,\dots,s_k)\bigr)=\frac1n\prod_{j=1}^{n-1}\Bigl(2k-\sum_{i=1}^k(\omega_n^{s_i j}+\omega_n^{-s_i j})\Bigr),
\]
where $\omega_n=e^{2\pi i/n}$. Specializing to $k=2$, $(s_1,s_2)=(1,q)$ and $n=p$ gives
\[
\tau(p,q)=\frac1p\prod_{j=1}^{p-1}\Bigl(4-2\cos\tfrac{2\pi j}{p}-2\cos\tfrac{2\pi q j}{p}\Bigr).
\]
Invariance under $q\mapsto -q$ and $q\mapsto q^{-1}$ follows since $j\mapsto qj$ permutes $\{1,\dots,p-1\}$ and $\cos$ is even. 
\end{proof}

\begin{theorem}\label{thm:square-shape}
Let $G=C_p(\pm 1,\pm q)$ with $\gcd(p,q)=1$, $p\ge3$. Then:
\begin{itemize}
\item[\emph{(a)}] If $p$ is odd, there exists an integer $A$ (depending on $p$, $q$) such that
\[
\tau(p,q)=p\,A^2.
\]
\item[\emph{(b)}] If $p$ is even, there exists an integer $B$ (depending on $p$, $q$) such that
\[
\tau(p,q)=\frac{\lambda_{p/2}}{p}\,B^2,
\qquad \lambda_{p/2}=6-2(-1)^q\in\{4,8\}.
\]
\end{itemize}
\end{theorem}

\begin{proof}
By Theorem~\ref{thm:product-and-parity},
\[
\tau(p,q)=\frac{1}{p}\prod_{j=1}^{p-1}\lambda_j,
\qquad 
\lambda_j=4-2\cos\frac{2\pi j}{p}-2\cos\frac{2\pi q j}{p}.
\]

\begin{itemize}
    \item[(a)] Suppose $p$ is odd. Hence, the product pairs up symmetrically:
\[
\prod_{j=1}^{p-1}\lambda_j=\left(\prod_{j=1}^{(p-1)/2}\lambda_j\right)^2,
\]
hence
\[
\tau(p,q)=\frac{1}{p}\left(\prod_{j=1}^{(p-1)/2}\lambda_j\right)^2.
\]
As $\tau(p,q)\in\mathbb{Z}$, the product must be divisible by $p$, say $\prod_{j=1}^{(p-1)/2}\lambda_j=p\,A$, giving $\tau(p,q)=p\,A^2$.

\item[(b)] Suppose $p$ is even. All eigenvalues pair except the unpaired one at $j=p/2$:
\[
\prod_{j=1}^{p-1}\lambda_j=\lambda_{p/2}\left(\prod_{j=1}^{(p/2)-1}\lambda_j\right)^2.
\]
Since,
\[
\lambda_{p/2}=4-2\cos\pi-2\cos(q\pi)=6-2(-1)^q\in\{4,8\}.
\]
We get,
\[
\tau(p,q)=\frac{\lambda_{p/2}}{p}\,B^2,
\qquad B=\prod_{j=1}^{(p/2)-1}\lambda_j \in\mathbb{Z}.
\]  
\end{itemize}
\end{proof}

\begin{remark}
In particular, for $q\equiv\pm1\pmod p$ one recovers $\tau(p,1)=p\,2^{p-1}$, consistent with Theorem~\ref{thm:square-shape}. Moreover, when $p$ is even the explicit expression $\lambda_{p/2}=6-2(-1)^q$ shows that the parity of $q$ directly affects the value of $\tau(p,q)$.
\end{remark}

The preceding results show that $T(G;x,y)$ already recovers $p$ and that $T(G;1,1)=\tau(p,q)$ depends only on the orbit $\{\pm q^{\pm1}\}$. Our computations up to $p\le 400$ suggest that the full two–variable polynomial carries enough additional structure to separate all homeomorphism classes.

\begin{conjecture}\label{conj:tau-injective} For every prime $p\ge3$, the map $q\mapsto \tau(p,q)$ is constant precisely on the orbits $\{\pm q^{\pm1}\}$. Equivalently, $\tau(p,q)$ distinguishes the homeomorphism classes of lens spaces $L(p,q)$. \end{conjecture}

{Since the Bollobás--Riordan polynomial specializes to the Tutte polynomial, it encodes strictly richer embedding information. In our data, rigid patterns in the exponents across the $z$–layers of such polynomials already distinguish the orbits $\{\pm q^{\pm1}\}$ and strongly point to full completeness.}

{\begin{conjecture}\label{conj:BR-complete}
The Bollobás--Riordan polynomial of the Heegaard graph of lens spaces is a complete invariant for lens spaces.
\end{conjecture}}

\rev{These conjectures bear directly on the possibility of proving the classification of lens spaces from the polynomials themselves. The parameter $p$ is read off the polynomial explicitly. The parameter $q$, on the other hand, can only be seen up to its orbit $\{\pm q^{\pm1}\}$: the symmetries realizing the equivalence $q\sim\pm q^{\pm1}$ become symmetries of the Heegaard ribbon graph, so the polynomial, being an invariant, can depend on $q$ only through this orbit. The subtler question is whether \emph{distinct} orbits are always separated by the polynomial; an affirmative answer would, conversely, yield a polynomial proof of the classification. We leave this as an open problem.}

\paragraph{The Penrose polynomial}

Building on Definition~\ref{df:pen} and Theorem~\ref{teo:val}, in this subsection we specialize to Heegaard graphs of lens spaces. We show that from the Penrose polynomial one can recover the value of $p$, and obtain partial information about the value of $q$. Some integer evaluations of Penrose polynomials of embedded graphs have been studied. In \cite{ELLISMONAGHAN2013424}, the number of vertices of an embedded checkerboard-colourable graph is \rev{obtained} from the evaluation of its Penrose polynomial at $2$; we extend the same result for embedded $4$-regular graphs.

\begin{theorem} \label{teo:ver}
Let $G$  be a cellularly embedded graph in an orientable closed surface, then we have: 
\begin{itemize}
\item[(i)] $P(G;1)=0$.
\item[(ii)] If $G$ has $p$ vertices and each vertex has degree $4$, then $P(G;2)=2^{p}$. 
\end{itemize}
\end{theorem}

\begin{figure}
  \centering
    \includegraphics{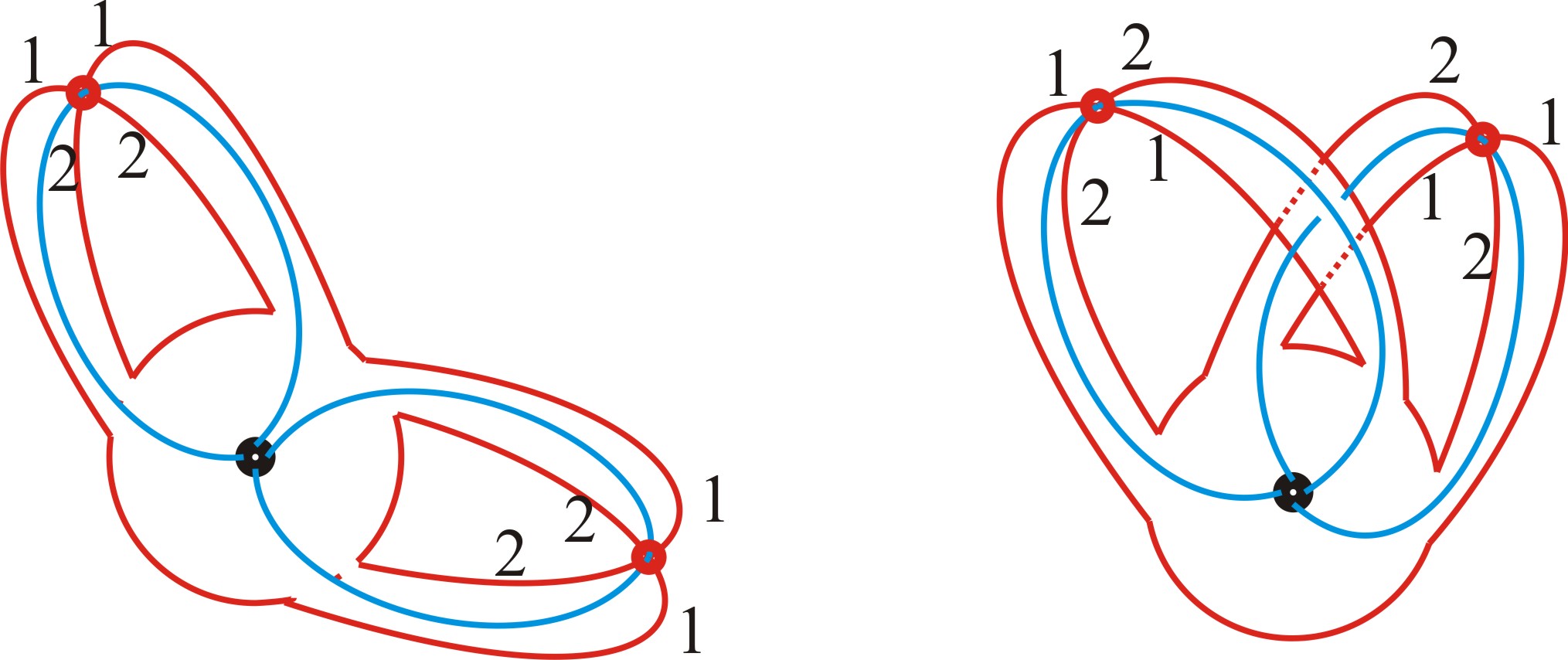}
  \caption{Admissible $2$-valuations}
  \label{fig:f8}
\end{figure}

\begin{figure}
  \centering
    \includegraphics{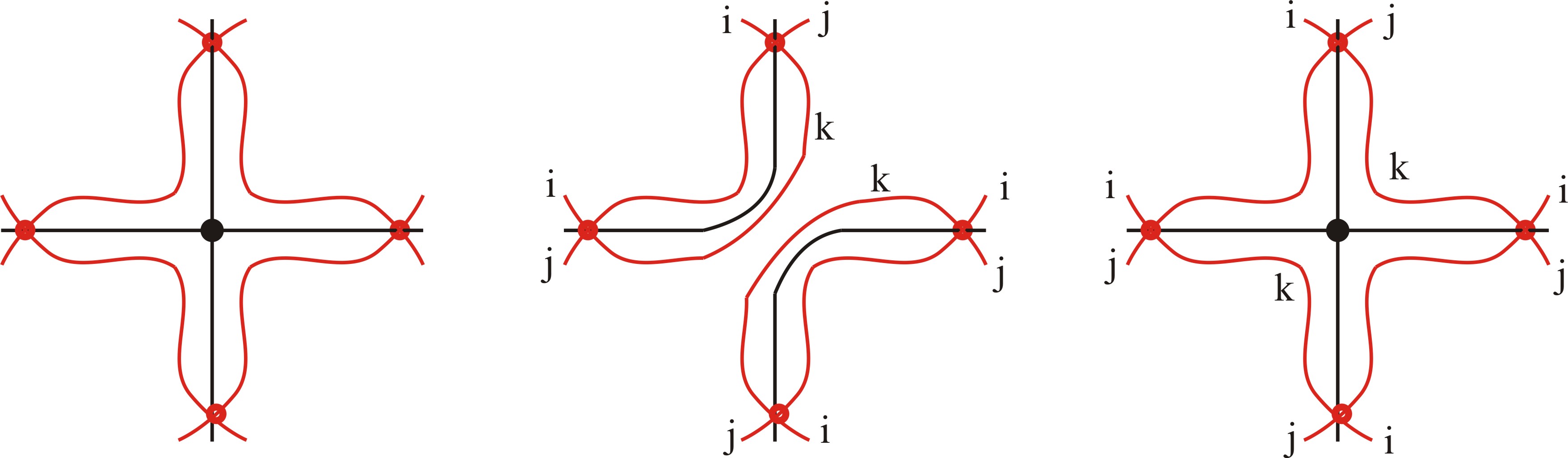}
  \caption{Admissible $2$-valuations of a reduced graph}
  \label{fig:f7}
\end{figure}

\begin{proof}
Let $G_m$ be the medial graph of $G$.
\begin{itemize}
\item[(i)] Order the edges of $G$, say $v_1,v_2,\ldots, v_r$. Every Penrose state $s$ of $G_m$ can be represented by a vector $(x_1,x_2,\ldots, x_r)$, where a coordinate $x_i\in \{0,1\}$ is $0$ if the smoothing of the state $s$ corresponding to the edge $v_i$ is a white smoothing, and it is $1$ if the smoothing is of crossing type. Let $s_1$ and $s_2$ be two Penrose states with identical vectors except for the last coordinate $x_r$, say $s_1$ has the last coordinate equal $1$, while the last coordinate of $s_2$ is $0$. The contributions of $s_1$ and $s_2$ to the Penrose polynomial have opposite signs, and in the evaluation $P(G;1)$, they cancel each other. Since every Penrose state can be paired with the state differing in the last coordinate, we have that $P(G;1)=0$.   
\item[(ii)] We proceed by induction on the number of vertices of $G$.  If $G$ has one vertex, since $G$ is embedded in an orientable surface and this vertex has degree $4$, the configuration of the medial graph $G_m$ corresponds to one of the two configurations in Fig.~\ref{fig:f8}. It \rev{is} easy to check that each one of these configurations has exactly  two admissible $2$-valuations: the one shown in Fig.~\ref{fig:f8} and reversing the two colors. We conclude from  Theorem \ref{teo:val} that $P(G;2)=2$.  \\
Now suppose that the evaluation at $2$ of  the Penrose polynomial of a $4$-regular graph with $m$ vertices  is equal to $2^{m}$. Let $G$ be a $4$-regular graph with $m+1$ vertices which is cellularly embedded in an orientable closed surface $F$. Let $v$ be a vertex of $G$. Construct a graph $\hat{G}$ from $G$ after removing the vertex $v$ and smoothing as shown in the middle diagram of Fig.~\ref{fig:f7}. Note that $\hat{G}$ is embedded  in $F$ and has $m$ vertices. We can also obtain the medial graph $\hat{G}_m$ from $G_m$ following the smoothing and considering each pair of connected vertices as a single vertex as shown in Fig.~\ref{fig:f7}. Given an admissible $2$-valuation of $\hat{G}_m$, it only remains to assign colors to the four arcs surrounding $v$. Note that the colors of the arcs labeled with $k$ in Fig.~\ref{fig:f7} should be the same in order to induce an admissible $2$-valuation of $G_m$. Since we have two possible choices for this selection, we conclude that $P(G;2)=2\cdot P(\hat{G};2)= 2 \cdot 2^{m}=2^{m+1}$.      
\end{itemize}
\end{proof}

In the case of Heegaard graphs obtained from a Heegaard decomposition of a $3$-manifold, since they are $4$-regular and are embedded in the orientable Heegaard surface, we conclude from Theorem \ref{teo:ver}:  
\begin{corollary} If $G$ \rev{is} a Heegaard graph of a $3$-manifold, then the number of vertices of $G$ is \rev{$\log_{2}$}$(P(G;2))$.
\end{corollary}

In order to show that the graph polynomial obtained from the Penrose polynomial is a complete invariant of lens spaces, it remains to recover the $q$ parameter from the polynomial. We present some advances in this direction by characterizing all the Penrose polynomials for lens spaces $L(p,1)$.

\begin{theorem} \label{teo:uno}
Let $G$ be the Heegaard graph of lens space $L(p,q)$, $1\leq q < p$. Then, $P(G;\lambda)$ is monic and $deg(P(G;\lambda))> p$ if and only if $q=1$.

\begin{figure}
  \centering
    \includegraphics{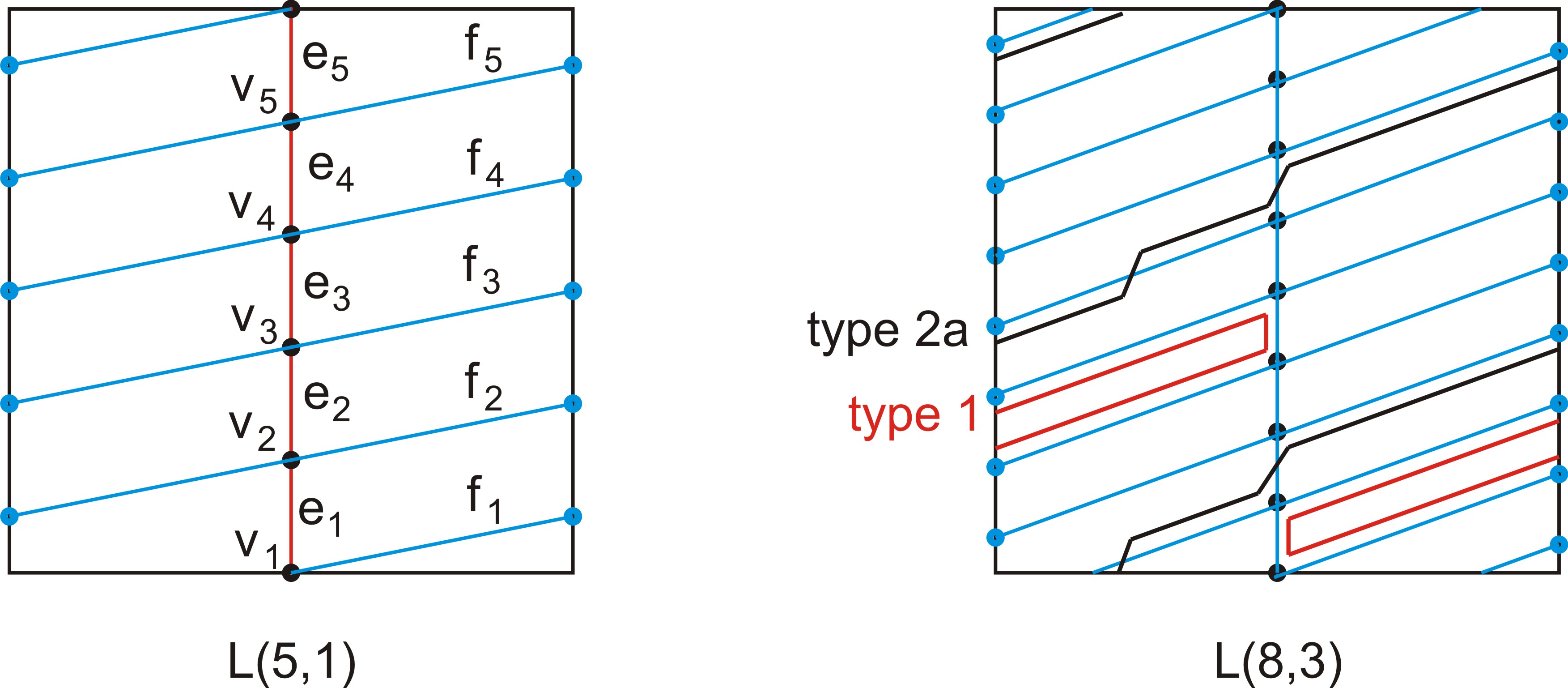}
  \caption{Left: enumeration of edges and vertices of $G$. Right: two types of length-$4$ boundaries in a Penrose state of $G_m$}
  \label{fig:f9}
\end{figure}

\begin{proof}

Let $D_1$ and $D_2$ be a system of decomposing disks of $L(p,q)$ to construct $G$ as described in Section~\ref{heegaardgraphs}. Label the edges $e_1,e_2, \ldots,e_p$ and $f_1,f_2, \ldots,f_p$ of $G$ in order as they appear around $\partial D_1$ and $\partial D_2$, respectively. Label the vertices of $G$, $v_1,v_2,\ldots,v_p$, such that $v_i$ is between $e_{i-1}$ and $e_{i}$. In the left-side picture in Fig.~\ref{fig:f9}  we show the enumeration of edges and vertices of the graph corresponding to $L(5,1)$.

In order to show that  $P(G;\lambda)$ is monic, we analyze the Penrose states of $G_m$ with maximum number of boundary components. Given $\gamma$ a boundary component of a Penrose state of $G_m$, note that $\gamma$ runs along some edges of $G$, alternating between the two families of arcs, namely, the path of $\gamma$ along the edges of $G$ could be $e_{i_1},f_{i_2},e_{i_3},\ldots, f_{i_r}$. In this case, we say that the length of $\gamma$ is $r$. Note that $r$ should be even, since $\gamma$ is a closed curve. An edge $e_i$ connects the vertices $v_{i}$ and $v_{i+1}$ (identifying $v_{p+1}$ with $v_1$), while the edge $f_i$ connects a vertex $v_j$ with $v_{j+q}$ (modulo $p$). In other words, the curve $\gamma$ starts at some vertex $v_j$ of $G$ and follows a progression in the indices of the vertices given by the sequence $\pm 1, \pm q, \ldots, \pm 1, \pm q$ (as many terms as $r$). This imposes the condition that $p$ should divide $\pm 1 \pm q \pm\cdots \pm 1 \pm q$, since $\gamma$ is a closed curve.

In the case that $G$ admits a Penrose state $s$ with a boundary component $\gamma$ of length $2$, then $p$ must divide $\pm 1\pm q$, and the only possibility is that $q=1$ (if  $q=p-1$, the lens spaces is equivalent to the one with $q=1$ according to Theorem \ref{teo:clas}). If a Penrose state of $G$ has a boundary of length $4$, then $p$ must divide $\pm 1\pm q \pm 1 \pm q$. The only possibilities for the signs are $1+q-1-q$, for any value $q$, or $1+q+1+q=p$ or $-1+q-1+q=p$. Note that some combinations such as $1+q-1+q =2q=p$ are not allowed since $p$ and $q$ are relatively prime. Two of these types of length-$4$ boundaries are illustrated in right-side diagram in Fig.~\ref{fig:f9}. Let us call a length-$4$ boundary components of a Penrose state of $G_m$ of type $1$, type $2a$ or type $2b$ if it corresponds to the solutions  $1+q-1-q=0$, $1+q+1+q=p$ or $-1+q-1+q=p$, respectively.     

First, consider the case $q=1$. Let $s$ be the Penrose state of $G_m$ with all the smoothings of crossing type. It is easy to check that this state has $p+1$ boundary components if $p$ is odd \rev{and} $p+2$ boundary components if $p$ is even. We proceed to prove that $s$ is the Penrose state of $G_m$ with the maximum number of boundary components and that it is unique with this property. Let $s'$ be another Penrose state and let $t$ be the number \rev{of} length-$2$ boundaries in $s'$. If $t=0$, then the length of every boundary of $s'$ is at least $4$ and the number of boundaries is at most $p$. In the case where  $p>4$, there are no length-$4$ boundaries of type $2a$ or $2b$, and there are at most $p-t-1$ boundaries of type $1$. Then the number of boundaries of $s'$ of length $2$ or $4$ is $t+(p-t-1-k)=p-k-1$, for some $k\geq 0$, and the sum of the lengths of this boundaries is $2t+4(p-t-1-k)=4p-2t-4k-4$. Since the sum of the lengths of all the boundary components in $s'$ is $4p$, then the sum of the lengths of the boundaries of length greater than $4$ is $2t+4k+4$. Increasing $k$ by $1$, decreases the number of boundaries of length $2$ or $4$ by $1$, but the number of boundaries of length greater than $4$ may not increase by $1$. This leads us to conclude that the maximum number of boundary components in $s'$ occurs when $k=0$ and $t$ is the maximum possible, namely, $t=p$. This later case occurs  when $s'=s$. In the case $p=4$, there are boundaries of $s'$ of type $2a$, and all Penrose states of $G_m$ can be analyzed to conclude that $s'$ is the state with a maximum number of boundaries and this number is $6$. 

Finally, consider the case $q>1$. Let $s$ be the Penrose state of $G_m$ with all the smoothings white. Note that this state has $p$ boundary components; as before, we will show that this is the unique state with the maximum number of boundary components. Let $s'$ be another Penrose state of $G_m$. If there are only length-$4$ boundaries in $s'$ of type $1$, then the number of boundaries in $s'$ is at most $p$, and it is $p$ only if $s'=s$. If $s'$ has $k\geq 1$ boundary components of type $2a$, then it does not have boundaries of type $2b$ (and vice versa), and it has at most $p-3-k$ boundaries of type $1$. Since $k$ is bounded above by $p/2$, we conclude that $s'$ has less boundary components than $s$. Then $s$ is the Penrose state with the largest number of boundary components and this number is $p$.  
\end{proof}
\end{theorem}

\subsection{A first look at the Poincaré sphere}\label{subsec:Poincare}

To compute the polynomial associated with a minimal Heegaard graph of a 3-manifold we need two things: to know a minimal genus Heegaard decomposition and a corresponding minimal Heegaard graph.

In general, finding minimal decompositions of Heegaard genus has been a research challenge. However, there are known results in this regard. For example, consider the Poincaré sphere $P^3$.

The Poincaré sphere has several characterizations, see \cite{KSch}, here we are interested in two of them: as a Seifert fibered space and as Heegaard genus two manifold.

Recall that $P^3$ is the Seifert fibered space with base space the 2-sphere with three exceptional fibers of orders $(2,1), (3,1), (5,1)$. In \cite{BoileauCollinsZieschang} and \cite{BoileauOtal} Boileau et al., classify the genus two Heegaard splittings of  Seifert fibered spaces with base space the 2-sphere and three exceptional fibers. In particular the uniqueness (up to isotopy and up to homeomorphism) of the genus two Heegaard splitting of $P^3$ is known. \rev{Hence $P^3$ has a unique genus-two Heegaard diagram up to the equivalence of Definition~\ref{def:diagram-equiv}.}

As a consequence of Theorem \ref{FicoTheorem1} the minimum number of vertices of a Heegaard graph for a homology sphere is at least twelve. \rev{The diagram of Fig.~\ref{fig:PoincareSphere} realizes this minimum and so yields a Heegaard graph of the Poincaré homology sphere $P^3$. To the best of our knowledge, it is the only such graph presently known. Evaluating the Penrose polynomial on it yields}
\begin{multline*}
z^{12} - 24z^{11} + 553z^{10} - 6186z^{9} + 42664z^{8} - 193904z^{7} \\
+ 595168z^{6} - 1238528z^{5} + 1718528z^{4} - 1518592z^{3} + 770816z^{2} - 170496z
\end{multline*}

\rev{This motivates the following conjecture.}

\rev{\begin{conjecture}\label{conj:poincare-unique}
The Poincaré homology sphere has a unique Heegaard graph up to the equivalence of Heegaard graphs; equivalently, the set $\mathcal{G}(P^3)$ consists of a single equivalence class.
\end{conjecture}}

\rev{If Conjecture~\ref{conj:poincare-unique} holds, the Penrose polynomial computed above is a genuine invariant of $P^3$. If, on the contrary, the Poincaré sphere admitted several inequivalent minimal Heegaard graphs, the collection of their polynomials ---taken as an unordered tuple, or suitably combined into a single derived polynomial--- would still provide a well-defined invariant.}

\rev{More generally, it is natural to ask whether the phenomenon behind Conjecture~\ref{conj:poincare-unique} is peculiar to the Poincaré sphere.}

\rev{\begin{question}\label{quest}
Does a Heegaard diagram realizing the minimal possible (global) number of vertices determine a unique Heegaard graph up to equivalence? Equivalently, do any two equivalent such diagrams have equivalent Heegaard graphs?
\end{question}}

\rev{A positive answer to this question, together with the uniqueness of the genus-two Heegaard diagram of the Poincaré sphere up to equivalence, would imply Conjecture~\ref{conj:poincare-unique}.}
\rev{More broadly, while the Heegaard genus of a $3$-manifold, and \emph{a fortiori} a minimal Heegaard graph realizing it, are hard to determine algorithmically, the present example shows that explicit Heegaard graphs can nevertheless be studied effectively in concrete cases. This suggests a possible route toward combinatorial invariants of $3$-manifolds that have so far remained unexplored.}

\section{Conclusion}

We have introduced a framework to derive polynomial invariants of closed, orientable $3$-manifolds from Heegaard diagrams by passing to their associated minimal Heegaard graphs. This approach yields a computable setting where classical graph polynomials such as the Tutte, Penrose, and Bollobás–Riordan polynomials can be evaluated to produce invariants of $3$-manifolds.

For lens spaces, we showed that the associated graphs are double–loop circulants, uniquely embedded in the Heegaard torus. Consequently, any graph polynomial on these graphs defines a genuine lens space invariant. Moreover, the Tutte and Bollobás–Riordan polynomials recover the parameter $p$ and computational evidence strongly suggests they also detect $q$, raising the possibility that these polynomials form complete invariants. The Penrose polynomial, in turn, provides independent information and characterizes the special case $q=1$.

\rev{Finally, for the Poincaré homology sphere ---whose genus-two splitting is unique and whose minimal Heegaard graphs have, by the argument of González-Acuña in the Appendix, at least twelve vertices--- we computed the Penrose polynomial of the only currently known minimal Heegaard graph, illustrating that these invariants are computable beyond lens spaces. We conjecture that this graph is unique, in which case any such polynomial is a genuine invariant of the manifold.} These results suggest a central open question: \emph{to what extent can polynomial invariants of Heegaard graphs provide a complete classification of $3$-manifolds?}

\bigskip

\backmatter

\bmhead{Acknowledgements}

We warmly thank our colleague Francisco González-Acuña for sharing with us the argument establishing the 12-vertex minimum, reproduced in the Appendix with his permission. \rev{We are also grateful to the anonymous referee, whose careful reading and remarks significantly improved the manuscript.}

FM-G was supported by Grant UNAM-PAPIIT IN113323.

JL-LM was supported by a SECIHTI postdoctoral fellowship and by grant CBF2023-2024-4059.

\clearpage 

\appendix
\subsection*{Appendix: The minimal number of vertices for Heegaard diagrams of Homology spheres by Francisco González Acuña.}\label{section:appendix}

\begin{theorem}\label{FicoTheorem1}
A minimal genus $g\geq 2$ Heegaard diagram of a homology sphere, other than $S^3$, has at least twelve vertices.
\begin{proof}
Let $M$ be a homology sphere and  $D=(S; v,w)$ be a genus $g$ Heegaard diagram \rev{with $m$ vertices}. Let $\langle x_1, x_2, \dots x_g : r_1, r_2, \dots, r_g\rangle$ be the corresponding presentation for $\pi_1(M)$, then $|r_1|+|r_2|+\dots |r_g|=m$. Recall that the fundamental group of a homology sphere is a perfect group. If $g=2$, then $m\geq 12$, otherwise Proposition \ref{prop:Fico1} implies that $\pi_1(M)$ is trivial, which is a contradiction. If $g\geq 3$ then $m\geq 12$, otherwise by Proposition \ref{prop:Fico2} $M$ would be of lower Heegaard genus, which is not possible. 
\end{proof}
\end{theorem} 

Here, $|\cdot|$ denotes the word length.

\begin{proposition}\label{prop:Fico1}
If $G=\langle x_1, x_2 : r_1, r_2 \rangle$ is a perfect group and the sum of relator lengths $|r_1|+|r_2|<12$, then $G$ is the trivial group.

\begin{proof}
We can assume that the relators $r_1$ and $r_2$ are cyclically reduced.

Since $G$ is a perfect group, its abelianization is trivial. This implies that the determinant of the exponent sum matrix $D$ must be $\pm 1$. The matrix $D$ is given by:
\[
D = \begin{pmatrix}
    s_1(r_1) & s_2(r_1) \\
    s_1(r_2) & s_2(r_2)
\end{pmatrix}
\]
where $s_i(r)$ is the exponent sum of the generator $x_i$ in the relator $r$.

It is a known result that if any relator $r_j$ has length $|r_j| < 5$ ($j=1,2$), the group $G$ must be abelian. Since $G$ is also perfect, it must be trivial.

We can therefore assume that $|r_1|=5$ and $|r_2| \in \{5, 6\}$.

For brevity, we will use the following notation: $1$ for $x_1$, $1'$ for $x_1^{-1}$, $2$ for $x_2$, and $2'$ for $x_2^{-1}$.

A relator cannot be a word of a single generator, such as $11111$ or $22222$, as this would imply $s_1(r_1)=5$ or $s_2(r_1)=5$, making $\det(D)$ divisible by 5. This contradicts $\det(D) = \pm 1$. If a conjugate of a relator or its inverse is $1$ times a power of $2$ or $2$ times a power of $1$, the group $G$ is cyclic and therefore trivial.

For the remaining cases, we can apply an automorphism $A$ of the free group $F_2 = \langle x_1, x_2 \rangle$ (of the form $A(x_i) = x_j^{\pm 1}$) to simplify $r_1$ without changing the length of $r_2$. After applying such an automorphism and possibly taking a conjugate or inverse of $A(r_1)$, its abelianization can be reduced to one of three forms. This leads to a case analysis:

\begin{itemize}
    \item[] \textbf{Case I:} The abelianization of $r_1$ is $11122$.
    \item[] \textbf{Case II:} The abelianization of $r_1$ is $122$.
    \item[] \textbf{Case III:} The abelianization of $r_1$ is $1$.
\end{itemize}

\textbf{Case I.}
We can assume $r_1 = 11122$. (The other element of length 5 with this abelianization, $11212$, leads to a cyclic group 
because $11212$ is primitive in $F$). The condition $\det(D) = \pm 1$ implies that the exponent sums of the second relator, $(s_1(r_2), s_2(r_2))$, must be either $(1, 1)$ or $(2, 1)$. The possible forms for $r_2$ are $212'12$, $2112'1'2$, $21'2'112$, $212'121'$, $212'1'21$, or $21'2'121$.
\begin{itemize}
    \item If $G = \langle 1,2: 11122, 212'12 \rangle = \langle 1,2: 12211, 12212' \rangle \Rightarrow 1=2' \Rightarrow G$ is cyclic.
    \item If $G = \langle 1,2: 11122, 2112'1'2 \rangle = \langle 1,2: 22111, 22112'1' \rangle \Rightarrow 1=2'1' \Rightarrow G$ is cyclic.
    \item If $G = \langle 1,2: 11122, 21'2'112 \rangle = \langle 1,2: 11221, 11221'2' \rangle \Rightarrow 1=1'2' \Rightarrow G$ is cyclic.
    \item For the remaining relators ($212'121'$, $212'1'21$, or $21'2'121$), a computation in GAP shows that the resulting group $\langle 1,2: r_1, r_2 \rangle$ is trivial.
\end{itemize}

\textbf{Case II.}
We can assume $r_1 = 21121'$. The condition on the determinant restricts $(s_1(r_2), s_2(r_2))$ to a few possibilities. The 8 possible words for $r_2$ are: $12'1'22$ (for sums $(0,1)$), $121'212$ (for sums $(1,3)$), $12122$ (for sums $(2,3)$), or $121'212'$, $121'2'12$, $12'1'212$, $1221'2'1$, $12'1'221$ (for sums $(1,1)$). A computation in GAP shows that all 8 of these group presentations yield the trivial group.

\textbf{Case III.}
We can assume $r_1 = 121'2'1$ (the commutator $[1,2]$ followed by $1$). For $r_2$, we can write it in the form $2a2'b2c$ where $a,b,c$ are powers of $1$ such that $|a|+|b|+|c| \in \{2,3\}$ with $|a|>0$ and $|b|>0$. This gives the following 20 possible words for $r_2$:
\begin{align*}
& 212'12, && 212'1'2, && 21'2'12, && 21'2'1'2,  \\
& 212'121, && 212'1'21, && 21'2'121, && 21'2'1'21, \\
& 212'121', && 212'1'21', && 21'2'121', && 21'2'1'21', \\
& 2112'12, && 2112'1'2, && 21'1'2'12, && 21'1'2'1'2, \\
& 212'112, && 21'2'112, && 212'1'1'2, && 21'2'1'1'2.
\end{align*}
A computation in GAP shows that all 20 of these corresponding group presentations yield the trivial group.

Since all cases where $|r_1|+|r_2|<12$ lead to a trivial group, the proposition is proven.
\end{proof}
\end{proposition}

\begin{proposition}\label{prop:Fico2}
    Let $M$ be a homology sphere with Heegaard genus $g\geq 3$. Let  $D=(S; v,w)$ be a Heegaard diagram associated to the Heegaard decomposition of genus $g$. And let $\langle x_1, x_2, \dots, x_g : r_1, r_2 \dots , r_g \rangle$ be the corresponding presentation of $\pi_1(M)$. Then $|r_1|+|r_2|+\dots |r_g| \geq 12$.

\begin{proof}
    Suppose that  $|r_1|+|r_2|+\dots |r_g| < 12$, then  $|r_i|<4$ for some $i=1,2, \dots, g$. After applying an automorphism $A$ of the free group $F$ generated by $x_1, x_2, \dots x_g$ such that $A(x_i)=x_j^{\pm1}$, we can assume that $r_1$ is $x_1x_2x_3$, $x_1x_2^2$, $x_1x_2$ or $x_1$. Notice that $r_1$ cannot be $x_1^2$ or $x_1^2$ since $M$ is a homology sphere. The four options for $r_1$ are primitive elements in $F$, we can apply another automorphism $A'$ of $F$ such that $A'(r_1)=x_1$. This implies that $M$ has a Heegaard diagram of genus $g-1$, which contradicts the minimality of the Heegaard genus. 
\end{proof}
\end{proposition}

\begin{remark}
    In a similar way, it can be proven that a homology sphere $M$ with Heegaard genus $g$ and with  $\langle x_1, x_2, \dots, x_g : r_1, r_2 \dots , r_g \rangle$ the corresponding presentation of $\pi_1(M)$. Then $|r_i|>4$ and $|r_1|+|r_2|+\dots |r_g| \geq 5g$
\end{remark}

\clearpage


\end{document}